\documentclass[12pt]{article}
\usepackage{amsmath}
\usepackage{amssymb}
\usepackage{amsthm}
\usepackage{url}
\newtheorem{thm}{Theorem}
\newtheorem{lem}{Lemma}


\begin{document}

\begin{center}
{\bf\Large McCulloch-Pitts brains\\\rule{0pt}{20pt} and pseudorandom functions}\\
\vspace{0.5cm} Va\v sek Chv\' atal \footnote{{\tt
    chvatal@cse.concordia.ca}}, Mark Goldsmith
\footnote{{\tt markgoldsmith@gmail.com}},  and Nan Yang
\footnote{{\tt nan.yang@me.com}}\\
Department of Computer Science and Software Engineering\\
 Concordia University, Montreal
\end{center}

\begin{center}
{\bf Abstract}
\end{center}
{\small In a pioneering classic, Warren McCulloch and Walter Pitts
 proposed a model of the central nervous system. Motivated by EEG
 recordings of normal brain activity, Chv\' atal and Goldsmith asked
 whether or not these dynamical systems can be engineered to produce trajectories which are
irregular, disorderly, apparently unpredictable. We show that they 
cannot build weak pseudorandom functions.
\vspace{0.5cm}

Electroencephalogram recordings of normal brain (or of an epileptic
brain well before a seizure) are usually irregular, disorderly, with
no apparent pattern: see, for instance,
\cite{Liu66,Leh01,DaS03,Ias04,Cha09,Oca09,ATE10}. Chv\' atal and
Goldsmith~\cite{ChvGol12} asked whether or not the McCulloch-Pitts
model of the brain can be engineered to exhibit similar behaviour.
The same question, although without its physiological interpretation, 
was also asked in~\cite{elyada2005can}. Let us begin by briefly describing 
the McCulloch-Pitts model.

A {\em linear threshold function\/} is a function $f:{\bf R}^n
\rightarrow \{0,1\}$ such that, for some real numbers
$w_1,\ldots,w_n$ and $\theta$,
\[
\textstyle{
f(x_1,\ldots,x_n)= H\left(\sum_{j=1}^n  w_jx_j - \theta\right)
}
\]
where $H$ is the Heaviside step function defined by $H(d)=1$ for all
nonnegative $d$ and $H(d)=0$ for all negative $d$.  Warren
McCulloch and Walter Pitts~\cite{McCPit43} proposed a model of the
central nervous system built from linear threshold functions. When
this system has $n$ neurons and no peripheral afferents, its
McCulloch-Pitts model is a mapping $\Phi:\{0,1\}^n\rightarrow
\{0,1\}^n$ defined by
\[
\Phi(x) = (f_1(x), \ldots, f_n(x))
\]
for some linear threshold functions $f_1,\ldots ,f_n$.  We will
refer to such mappings $\Phi$ as {\em McCulloch-Pitts dynamical systems}.

Chv\' atal and Goldsmith~\cite{ChvGol12} asked whether or not these
dynamical systems can produce trajectories which are irregular,
disorderly, apparently unpredictable in the sense of generating random
numbers. In making the meaning of their question precise, they took
the point of view of the practitioners, who mean by a random number
generator any deterministic algorithm that, given a short sequence of
numbers, called a {\em seed,\/} returns a longer sequence of numbers;
such a random number generator is considered to be good if it passes
statistical tests from some commonly agreed on battery. (This point of
view is expounded in~\cite[Chapter 3]{Knu14}.) 
	
In this note, we take the point of view of the theorists: we are going to prove that McCulloch-Pitts dynamical systems cannot produce trajectories which are irregular,
disorderly, apparently unpredictable in the sense of providing {\em weak pseudorandom functions.\/} These 
have been introduced in~\cite{naor1995synthesizers} and subsume {\em pseudorandom functions,\/} introduced in~\cite{GGM86} under the original name of 
`poly-random collections'. Roughly speaking, a weak pseudorandom function is a probability distribution on a set $F_n$ of functions 
from $\{0,1\}^n$ to $\{0,1\}^n$ with the following property: if $x^1, \ldots ,x^m$ are chosen independently and uniformly at random from 
$\{0,1\}^{n}$,
then no polynomial-time randomized algorithm can distinguish with a non-negligible probability between 
(i) a sequence $(x^1,f(x^1), \ldots ,x^m,f(x^m))$ where $f$ is chosen at random from $F_n$ and (ii) a sequence 
$(x^1,y^1,\ldots ,x^m,y^m)$ where
$y^1,\ldots ,y^m$ are chosen independently and 
  uniformly at random from $\{0,1\}^{n}$. 
	(Distinguishing between (i) and (ii) is a trivial matter when $f$ is known and that is why an unknown $f$ must be drawn from  
	a probability distribution on $F_n$.)
	Our result shows that weak pseudorandom functions cannot be built from McCulloch-Pitts dynamical systems:

\begin{thm}\label{thm.main}
There is a polynomial-time deterministic algorithm that, given 
a sequence $(x^1,y^1,\ldots , x^m,y^m)$ of $n$-bit vectors, 
returns either the message {\tt McCulloch-} {\tt Pitts}
or the message {\tt not McCulloch-Pitts} in such a way that
\begin{itemize}
\item[{\rm (i)}] if $y^1=\Phi(x^1),  \ldots
  ,y^m=\Phi(x^m)$ for some McCulloch-Pitts dynamical system
  $\Phi$, then the algorithm returns {\tt McCulloch-Pitts},
\item[{\rm (ii)}] if $x^1, \ldots ,x^m$ are chosen independently and
  uniformly at random from $\{0,1\}^{n}$, if $y^1, \ldots ,y^m$ are chosen independently and
  uniformly at random from $\{0,1\}^{n}$, and if $m\ge (2+\varepsilon)n$ for some positive constant
  $\varepsilon$, then the algorithm returns 
	{\tt not McCulloch-Pitts} with probability at least $1-e^{-\delta n}$, where
  $\delta$ is a positive constant depending only on $\varepsilon$.
\end{itemize}
\end{thm}

A {\em dichotomy} of a set $X$ is its partition into two disjoint
sets. Unlike Cover~\cite{Cov65}, for whom a dichotomy is an unordered
pair of sets, we view every dichotomy as an ordered pair of sets. A
dichotomy $(X^+, X^-)$ of a subset of $\mathbf{R}^n$ is {\em linearly
  separable} if there are numbers $y_1,\ldots ,y_{n+1}$ such that
\begin{equation}\label{lsd}
\begin{split}
\textstyle{\sum}_{j=1}^n x_j y_j > y_{n+1} &\;\;\text{ whenever }
(x_1,\ldots ,x_n)\in X^+,\\
\textstyle{\sum}_{j=1}^n x_j y_j < y_{n+1} &\;\;\text{ whenever }
(x_1,\ldots ,x_n)\in X^-.
\end{split}
\end{equation}
When $f$ is a function from $\{0,1\}^n$ to $\{0,1\}$ and $x^1, \ldots ,x^m$ are points in
 $\{0,1\}^{n}$, the dichotomy $(\{x^i: f(x_i)=0\}, \{x^i: f(x_i)=1\})$ is linearly separable if and only if $f$ is a threshold function. Our proof of Theorem~\ref{thm.main} evolves from the propositions that 
linearly separable dichotomies are easy to recognize and linearly separable dichotomies are rare:
\begin{lem}
\label{dichlemma}
Linearly separable dichotomies of $m$-point subsets of $\{0,1\}^n$ can 
be recognized in time polynomial in $m$ and $n$.
\end{lem}
\begin{lem}
\label{problemma}
For every positive $\varepsilon$ there is a positive $\gamma$ with the
following property: If $\;X$ is a finite subset of $\mathbf{R}^n$ such that
$\lvert X\rvert\ge (2+\varepsilon)n$, then a dichotomy chosen uniformly at random from all
dichotomies of $X$ is linearly separable with probability at most
$e^{-\gamma n}$.
\end{lem}

Following the seminal report~\cite{RB}, the subject of 
learning a hyperplane that separates, or at least nearly separates, the two parts of a 
dichotomy received much attention in the machine learning 
community. None of it is relevant to the following standard argument, 
implicit in the linear programming proof of Minkowski's Separating Hyperplane Theorem 
for convex polytopes~\cite{tucker1955linear}.

\bigskip\noindent{\itshape Proof of Lemma~\ref{dichlemma}. }
Deciding whether a prescribed dichotomy of an $m$-point subset of $\{0,1\}^n$ 
is linearly separable amounts to solving system \eqref{lsd} of $m$ strict
linear inequalities in variables $y_1,\ldots ,y_{n+1}$, where each coefficient 
$x_j$ is $0$ or $1$; the epoch-making result of
Khachiyan \cite{Kha79} guarantees that this can be done in time polynomial in $m$ and $n$.
\hfill$\Box$\medskip

\bigskip\noindent{\itshape Proof of Lemma~\ref{problemma}. }
Without loss of generality, we may assume that 
that $0<\varepsilon\le 1$. Let $m$ denote $\lvert X\rvert$ and
let $p$ denote the probability that a dichotomy chosen uniformly at random from
all dichotomies of $X$ is linearly separable. 

Of the  $2^m$ dichotomies of $X$, at most
$2 \sum_{i=0}^{n} \binom{m-1}{i}$ are linearly separable (this is at least implicit in \cite{Win66} and \cite{Cov65}), and so  
\[\textstyle{
p\le 2^{-m+1} \sum_{i=0}^{n} \binom{m-1}{i}\le 2^{-m+1} \sum_{i=0}^{n} \binom{m}{i}.
}
\]
Since $m\ge (2+\varepsilon)n$ and
$0<\varepsilon\le 1$, we have $n\le (0.5-\varepsilon/6)m$;  a special case of the
well-known bound on the tail of the binomial distribution (see, for
instance,~\cite[Theorem 1]{Hoe63}) guarantees that for every positive
$\alpha$ smaller than $0.5$ there is a positive $\beta$ such that
\[
\textstyle{
\sum_{i\le (0.5-\alpha)m} \binom{m}{i}
\;\le\;
2^me^{-\beta m};
}
\]
setting $\alpha=\varepsilon/6$, we conclude that $p\le 2e^{-\beta m}$, which
proves the lemma.
\hfill$\Box$\medskip

An alternative proof of Lemma~\ref{problemma}, proposed by one of the reviewers, relies 
on the Sauer-Shelah Lemma (\cite{sauer1972density}, \cite{shelah1972combinatorial}):
{\em If a family of subsets of an $m$-point set has Vapnik-Chervonenkis dimension $d$,
then it includes at most $\sum_{i=0}^{d} \binom{m}{i}$ sets.\/} Its other ingredient is 
the following corollary of Radon's theorem \cite{Rad21}: {\em If $\mathcal{H}$ is a family 
of half-spaces in $\mathbf{R}^n$ and if $X$ is a finite subset of $\mathbf{R}^n$, then family 
$\{X\cap Y: Y\in \mathcal{H}\}$ has Vapnik-Chervonenkis dimension at most $n+1$.\/} Putting 
the two together, we conclude that 
$X$ has at most $2 \sum_{i=0}^{n+1} \binom{m}{i}$ linearly separable dichotomies.  
This upper bound, although weaker than our $2 \sum_{i=0}^{n} \binom{m-1}{i}$, also 
yields the lemma's conclusion.

\bigskip\noindent{\itshape Proof of Theorem~\ref{thm.main}. }
The algorithm goes as follows: Let $\alpha^i$ denote the first bit of $y^i$ and
define 
\begin{align*}
X^+ &=\{x^i: 1\le i\le m,\; \alpha^i=1\},\\
X^- &=\{x^i: 1\le i\le m,\; \alpha^i=0\}.
\end{align*}
If this dichotomy is linearly separable, then return {\tt
  McCulloch-Pitts}; else return {\tt not McCulloch-Pitts}.

Lemma~\ref{dichlemma} guarantees that the algorithm can be implemented to run in polynomial time.

To prove (i), assume that $y^1=\Phi(x^1),  \ldots
  ,y^m=\Phi(x^m)$ for some McCulloch-Pitts dynamical system
  $\Phi:\{0,1\}^n\rightarrow \{0,1\}^n$ defined
by $\Phi(x) = (f_1(x), \ldots, f_n(x))$. Now 
$\alpha^i=f_1(x^i)$ for all $i=1,\ldots ,m$, which means that $f_1$
takes value $1$ on all points of $X^+$ and value $0$ on all points of
$X^-$; since $f_1$ is a threshold function, the dichotomy $(X^+, X^-)$
is linearly separable, and so the algorithm returns {\tt
  McCulloch-Pitts}.

To prove (ii), assume that $x^1, \ldots ,x^m$ are chosen independently and
  uniformly at random from $\{0,1\}^{n}$, that $y^1, \ldots ,y^m$ are chosen independently and
  uniformly at random from $\{0,1\}^{n}$, and that $m\ge (2+\varepsilon)n$ for some positive constant
  $\varepsilon$. 	
Since the probability that the algorithm returns {\tt not McCulloch-Pitts}
increases as $m$ increases, we may replace the assumption that $m\ge
(2+\varepsilon)n$ by the assumption that $m= \lceil(2+\varepsilon)n\rceil$.
Write $X=X^+\cup X^-$. Since $x^1, \ldots ,x^m$ are chosen independently and
  uniformly from $\{0,1\}^{n}$, they are pairwise distinct with probability
$2^n(2^n-1)\cdots(2^n-m+1)/2^{nm}$;
	since
	\[
\frac{2^n(2^n-1)\cdots(2^n-m+1)}{2^{nm}}\ge \left(\frac{2^n-m}{2^n}\right)^m
=\left(1-\frac{m}{2^n}\right)^m\ge 1-\frac{m^2}{2^n},
\]
this probability is at least $1-5n^22^{-n}$.  When
$\lvert X\rvert=m$, the assumption that $y^1, \ldots ,y^m$ are chosen independently and
  uniformly from $\{0,1\}^{n}$ implies that the dichotomy $(X^+, X^-)$ of $X$ is chosen
uniformly from all dichotomies of $X$, in which case
Lemma~\ref{problemma} guarantees that $(X^+, X^-)$ is linearly
separable with probability at most $e^{-\gamma n}$ for
some positive constant $\gamma$ depending only on $\varepsilon$.  We conclude that
the algorithm returns {\tt not McCulloch-Pitts} with probability at least
$1-5n^22^{-n}-e^{-\gamma n}$, which is at least $1-e^{-\delta n}$ for
some positive constant $\delta$ depending only on $\varepsilon$.
\hfill$\Box$\medskip

There is an obvious refinement of the algorithm used in the proof of Theorem~\ref{thm.main}:
with $y^i_j$ standing for the $j$-th bit of $y^i$, test each of the $n$ dichotomies
\[
\left(\{x^i: 1\le i\le m,\; y^i_j=1\}, \{x^i: 1\le i\le m,\; y^i_j=0\}\right)
\phantom{xxx} (j=1,\ldots ,n)
\]
and return {\tt  McCulloch-Pitts} if and only if all $n$ of them are linearly separable. 
In the context of distinguishing McCulloch-Pitts functions from truly random functions, the extra work required in this refinement is pointless. The probability of returning {\tt  McCulloch-Pitts} when $y^1, \ldots ,y^m$ are chosen independently and uniformly at random from $\{0,1\}^{n}$ is at most $e^{-\delta n}$ in the original version and that is good enough; reducing it further to $e^{-\delta n^2}$ in the refinement is nice, but unnecessary. In addition, the assumption $m\ge (2+\varepsilon)n$ cannot be significantly relaxed even in the refinement: it is at least implicit in \cite{Win66} and \cite{Cov65} that a dichotomy chosen uniformly at random from
all dichotomies of a set of fewer than $(2-\varepsilon)n$ points in $\mathbf{R}^n$ is linearly separable with probability at least $1-e^{-\delta n}$.

Theorem~\ref{thm.main} implies that certain simple devices (namely, McCulloch-Pitts dynamical 
systems) cannot generate pseudorandomness. 
In the opposite direction, it has been proved that certain simple devices can generate pseudorandomness: examples can be found in~\cite{naor1999distributed}, \cite{krause2001pseudorandom}, \cite{nielsen2002threshold}, \cite{naor2004number}, \cite{applebaum2010cryptography}.

The question whether McCulloch-Pitts networks can produce trajectories which are
irregular, disorderly, apparently unpredictable remains open: all
depends on the interpretation of the terms ``irregular, disorderly,
apparently unpredictable''. When clinical neurologists visually
inspect an electroencephalogram, their vague criteria for declaring it
random-like are a far cry from the distinguishers that cryptographers use
to separate deterministic sequences from random sequences. As Avi
Wigderson~\cite[page 6]{Wig09} put it,
\begin{quote}
``Randomness is in the eye of the beholder, or more precisely, in its
computational capabilities ... a phenomenon (be it natural or
artificial) is deemed ``random enough,'' or pseudorandom, if the class
of observers/applications we care about cannot distinguish it from
random!''
\end{quote}
Many examples of generators that appear random
to observers with restricted computational powers are known. In
particular, pseudorandom generators for polynomial size constant depth
circuits have been constructed in~\cite{AjtWig85};
later, this work was greatly simplified and improved in~\cite{Nis91}.  O'Connor~\cite{OCo88} proved that an infinite
binary sequence appears random to all finite-state machines if and
only if it is $\infty$-distributed. Pseudorandom generators for
space-bounded computation have been constructed in~\cite{Nis92}. 
It is conceivable that McCulloch-Pitts dynamical systems could fool
neurologists into finding their trajectories unpredictable just as
they find normal electroencephalograms unpredictable. Proving this in
a formal setting with a suitable definition of `neurologists' is an
interesting challenge. 

A variation on our theme comes from the idea that 
in a brain of $n$ neurons, only $m$ neurons may be visible to the observer and
the remaining $n-m$ are hidden from view. Formally, given positive integers $m,n$ such that $m\le n$
and given a McCulloch-Pitts dynamical system $\Phi: \{0,1\}^{n}
\rightarrow \{0,1\}^{n}$, we may consider the mapping 
${\Phi}_m: \{0,1\}^{n}\rightarrow \{0,1\}^m$ 
such that ${\Phi}_m(x)$ is the $m$-bit prefix of $\Phi(x)$.
Can such mappings provide pseudorandomness? Our Theorem~\ref{thm.main} shows that the
answer is negative when $m=n$; one of the reviewers argued that, under the usual assumption that 
one-way functions exist, the answer is close to affirmative when $m=1$. Here is the argument: Every one-way function $f$ (as every boolean function) can be computed by a threshold circuit~\cite[Chapter 7]{parberry1994circuit}. When this circuit has $n$ gates and depth $d$, it can be imbedded in a McCulloch-Pitts dynamical system $\Phi: \{0,1\}^{n}\rightarrow \{0,1\}^{n}$, where the results of its computation show up with the time delay of $d$ units. Now $f$ is represented in the $d$-fold iteration of $\Phi$ 
and there is an appropriate projection $\pi$, the {\em hard-core bit\/} \cite{goldreich1989hard}, such
that the sequence $\pi(x)$, $\pi(f(x))$, $\pi(f(f(x)))$, \ldots is pseudorandom. 

Statistical properties of ${\Phi}_1$ have been studied in \cite[Section 4.2]{goldsmith15}. 
For instance, there is a McCulloch-Pitts dynamical system $\Phi: \{0,1\}^{37}
\rightarrow \{0,1\}^{37}$ such that the restriction of the trajectory of $\Phi$ on the first bit 
passes all ten statistical tests of the battery {\tt SmallCrush} implemented in
  the software library TestU01 of \cite{EcuSim07,EcuSim09}.

\begin{center}
{\bf Acknowledgments}
\end{center}
This research was undertaken, in part, thanks to funding from the
Canada Research Chairs program and from the Natural Sciences and
Engineering Research Council of Canada. We are grateful to P\' eter
G\' acs for helpful comments on a draft of this note and to Avi
Wigderson for telling us about Nisan's papers~\cite{Nis91, Nis92}. We 
also thank the two anonymous reviewers for their thoughtful comments that 
made us improve the presentation considerably.

\end{document}